# SOME RESULTS ON $2^{N-M}$ DESIGNS OF RESOLUTION IV WITH (WEAK) MINIMUM ABERRATION

### By Hegang H. Chen and Ching-Shui Cheng[1]

*University of Maryland and University of California, Berkeley*

It is known that all resolution IV regular $2^{n-m}$ designs of run size $N = 2^{n-m}$ where $5N/16 < n < N/2$ must be projections of the maximal even design with $N/2$ factors and, therefore, are even designs. This paper derives a general and explicit relationship between the wordlength pattern of any even $2^{n-m}$ design and that of its complement in the maximal even design. Using these identities, we identify some (weak) minimum aberration $2^{n-m}$ designs of resolution IV and the structures of their complementary designs. Based on these results, several families of minimum aberration $2^{n-m}$ designs of resolution IV are constructed.

**1. Introduction.** Fractional factorial designs, especially those with two-level factors, have a long history of successful use in scientific investigations and industrial experiments. A $2^{-m}$th fraction of a $2^n$ factorial design, consisting of $2^{n-m}$ distinct combinations, is referred to as a $2^{n-m}$ fractional factorial design. Such a design is called regular if it can be constructed by using a defining relation. How to choose a good fractional factorial design is an important issue. Minimum aberration (MA), introduced by Fries and Hunter ([1980](#)), has become the most popular criterion for selecting fractional factorial designs.

For a regular fractional factorial design $D$, each interaction that appears in the defining relation is called a defining word, and the resolution of the design is defined as the length of the shortest defining word. For each positive integer $i$, let $A_i(D)$ be the number of defining words of length $i$. Then, the resolution is equal to the smallest $i$ such that $A_i(D) > 0$. The vector $W(D) = (A_1(D), A_2(D), \ldots, A_n(D))$ is called the wordlength pattern of $D$.

Received April 2007; revised September 2008.

[1]Supported by NSF Grant DMS-05-05556.

*AMS 2000 subject classification.* 62K15.

*Key words and phrases.* Even design, minimum aberration, regular fractional factorial design, resolution, wordlength pattern.

The minimum aberration criterion chooses a design by sequentially minimizing $A_1(D), A_2(D), A_3(D), \ldots$. A $2^{n-m}$ design with maximum resolution $R_{\max}$ is said to have weak minimum aberration if it has the minimum number of words of length $R_{\max}$. A two-level regular design is called even if all its defining words are of even length. Throughout this paper, we denote the run size $2^{n-m}$ by $N$. It is well known that there is a unique (up to isomorphism) resolution IV design with $n = N/2$. This design is even and has the property that every even design is its projection onto a certain subset of factors [see Chen and Hedayat (1998b)]. For this reason, we call the resolution IV design with $n = N/2$ the maximal even design.

Chen (1998) studied the connection between wordlength patterns and projections of $2^{n-m}$ designs and showed that minimum aberration designs have good projection properties. Meanwhile, Cheng, Steinberg and Sun (1999) provided some insight into minimum aberration, and justified this criterion by demonstrating that it is a good surrogate for some model-robustness criteria. In recent years, there has been much progress in the construction of MA fractional factorial designs. Each regular $2^{n-m}$ design of resolution III or higher can be constructed by choosing $n$ factors from the saturated regular design of resolution III, which has $N - 1$ factors. The $N - 1 - n$ factors that are not chosen form another design called the complementary design. Chen and Hedayat (1996) and Tang and Wu (1996) established identities relating the wordlength pattern of a design to that of its complementary design. These identities can be used to help construct MA designs by choosing appropriate complementary designs. This approach is very useful when the complementary designs have a small number of factors. Chen and Hedayat (1998a) have identified all minimum aberration $2^{n-m}$ designs whose complementary designs have fewer than 64 factors.

When $n \leq N/2$, the complementary design theory developed by Chen and Hedayat (1996) and Tang and Wu (1996) is not useful, since, in this case, the complementary design has more factors than the original design. Thus, an alternative theory is needed. Maximal even designs are known to be MA for $n = N/2$. When $n < N/2$, MA designs must be of resolution IV or higher. Recent results in the literature of projective geometry [Bruen, Haddad and Wehlau (1998) and Bruen and Wehlau (1999)], as summarized in Chen and Cheng (2004) and Chen and Cheng (2006), showed that all regular $2^{n-m}$ designs of resolution IV with $5N/16 < n < N/2$ are projections of the maximal even design. Therefore, for $5N/16 < n < N/2$, one can also construct an MA design by appropriately deleting a subset of factors from the maximal even design. Thus, it is useful to develop a corresponding complementary design theory for maximal even designs, as Chen and Hedayat (1996) and Tang and Wu (1996) did for saturated designs of resolution III. Butler (2003) employed this idea to classify MA $2^{n-m}$ designs of resolution IV.

The objective of this article is to establish an explicit complementary design theory for maximal even designs and to further investigate minimum aberration $2^{n-m}$ designs with $5N/16 < n < N/2$. In Section 2, after introducing some technical tools, we derive combinatorial identities that relate the wordlength pattern of a $2^{n-m}$ even design and that of its complement in the maximal even design. These identities have explicit forms, so that the wordlength pattern of an even design can be readily calculated from that of its complement in the maximal even design. Using these results, we identify some weak minimum aberration designs of resolution IV and the structures of their complementary designs in Section 3. A lower bound on the minimum number of words of length four is derived in Section 4. Finally, several families of minimum aberration $2^{n-m}$ designs of resolution IV are constructed in Section 5.

**2. Some technical tools and results.**   Let $D$ be a $2^{n-m}$ regular fractional factorial design of resolution III or higher. Following the notation in Chen and Hedayat (1996), the treatment combinations in $D$ can be represented as row vectors as follows:

$$D = \{\mathbf{x} : \mathbf{x} = \mathbf{v}\mathbf{B}_n, \mathbf{v} \in V_{n-m}\}, \tag{1}$$

where $\mathbf{B}_n$ is an $(n-m) \times n$ matrix of rank $n-m$ over the finite field $\mathrm{GF}(2)$ and $V_{n-m}$ is the $(n-m)$-dimensional vector space over $\mathrm{GF}(2)$. The matrix $\mathbf{B}_n$ is called the *factor representation* of $D$. Let $k = n - m$ and $N = 2^k$. Then, each column of $\mathbf{B}_n$ can be identified with a point of $\mathrm{PG}(k-1, 2)$, where $\mathrm{PG}(k-1, 2)$ is the projective geometry of dimension $k-1$ over $\mathrm{GF}(2)$. So, a regular fractional factorial design as in (1) is determined by a set of $n$ points of $\mathrm{PG}(k-1, 2)$, say $T = \{\mathbf{a}_1, \ldots, \mathbf{a}_n\}$. When $n = N - 1$, we obtain the saturated regular design of resolution III by choosing all the $N - 1$ points of $\mathrm{PG}(k-1, 2)$. For $n < N - 1$, $T$ can also be obtained by deleting $N - 1 - n$ points from $\mathrm{PG}(k-1, 2)$. Without loss of generality, we can represent the $N - 1$ points of $\mathrm{PG}(k-1, 2)$ as

$$\underbrace{\mathbf{a}_1, \ldots, \mathbf{a}_n}_{T}, \underbrace{\mathbf{a}_{n+1}, \ldots, \mathbf{a}_{N-1}}_{\overline{T}}, \tag{2}$$

where $\overline{T} = \{\mathbf{a}_{n+1}, \ldots, \mathbf{a}_{N-1}\}$ consists of the points in $\mathrm{PG}(k-1, 2) \setminus T$. Let $D$ and $\overline{D}$ be the two fractional factorial designs corresponding to $T$ and $\overline{T}$, respectively. We call $\overline{D}$ the complementary design of $D$ in the saturated regular design of resolution III.

In the above geometric representation, it can be shown that a maximal even design corresponds to the complement of a $(k-2)$-dimensional projective geometry. Specifically, let $F$ be a set of $N/2 - 1$ points that is itself a $\mathrm{PG}(k-2, 2)$. Then, $\overline{F}$ gives the factor representation of a maximal even design, and $\overline{F} = \{\mathbf{a}, \mathbf{a} + F\}$ for any $\mathbf{a} \in \overline{F}$. We also point out that there is a one-to-one correspondence between the points in $\overline{F}$ and those in a $(k-1)$-dimensional Euclidean geometry over GF(2), an EG$(k-1, 2)$. If $D$ is an even design, then, since it must be a projection of the maximal even design, its corresponding $T$ can be considered as a subset of $\overline{F}$. In this case, PG$(k-1, 2)$ in (2) can be further partitioned into three parts,

$$(3) \quad \mathrm{PG}(k-1, 2) = \{\underbrace{\mathbf{a}_1, \ldots, \mathbf{a}_n}_{T}, \underbrace{\mathbf{a}_{n+1}, \ldots, \mathbf{a}_{n+N/2-1}}_{F}, \underbrace{\mathbf{a}_{n+N/2}, \ldots, \mathbf{a}_{N-1}}_{\widetilde{T}}\}.$$

The vectors in the $k$-dimensional linear space generated by the rows of (3) can be displayed as an $N \times (N-1)$ matrix

$$(4) \quad \{\underbrace{\mathbf{c}_1, \ldots, \mathbf{c}_n}_{D}, \underbrace{\mathbf{c}_{n+1}, \ldots, \mathbf{c}_{n+N/2-1}}_{D^F}, \underbrace{\mathbf{c}_{n+N/2}, \ldots, \mathbf{c}_{N-1}}_{\widetilde{D}}\}.$$

We call $\widetilde{D}$ the complementary design of $D$ in the maximal even design, while $\overline{D} = D^F \cup \widetilde{D}$ is its complement in the saturated regular design of resolution III. Note that $\widetilde{D}$ is also an even design and that both $D$ and $\widetilde{D}$ are of resolution IV or higher.

As pointed out earlier, when $n < N/2$, $\overline{D}$ has more factors than $D$; therefore, in this case, the results of Chen and Hedayat (1996) and Tang and Wu (1996) do not help alleviate the complexity of identifying MA $2^{n-m}$ designs. On the other hand, when $5N/16 < n < N/2$, since the MA designs must be even, in view of (4), one can construct $D$ via $\widetilde{D}$. Butler (2003) studied the relationship between the wordlength patterns of $D$ and $\widetilde{D}$ in terms of the moments of their design matrices. However, no explicit identities linking these wordlength patterns are available.

Employing techniques of Chen and Cheng (1999), we now derive explicit combinatorial identities that govern the relationship between the wordlength pattern of a regular $2^{n-m}$ even design $D$ and that of its complement $\widetilde{D}$ in the maximal even design. Such identities play an important role in classifying (weak) MA designs in Section 3.

For any positive integer $n$, the Krawtchouk polynomial is defined by

$$P_i(j; n) = \sum_{s=0}^{i} (-1)^s \binom{j}{s} \binom{n-j}{i-s}, \qquad i = 0, 1, 2, \ldots.$$

As in Chen and Cheng (1999), we define

$$(5) \quad \alpha_r(s) = \begin{cases} \dfrac{1}{2^r}\left(\dbinom{2^r-2}{s} + P_s(2^{r-1}; 2^r-2)(2^{r-1}-1) \right. \\ \qquad\qquad\qquad \left. - P_s(2^{r-1}-1; 2^r-2)2^{r-1}\right), \\ \qquad \text{for } s = 2, \ldots, 2^r - 2, \\ 1, \qquad \text{for } s = 1, \end{cases}$$

(6) $$\gamma_r(s) = \frac{1}{2^r}\left(\binom{2^r-1}{s} + P_s(2^{r-1}; 2^r-1)(2^r-1)\right)$$

$$\text{for } s = 0, \ldots, 2^r - 1.$$

Since $D$ and $\widetilde{D}$ in (4) are even designs, $A_i(D) = A_i(\widetilde{D}) = 0$ for all odd $i$. Using Theorem 1 of Chen and Cheng (1999), we can prove the following theorem, which links the wordlength patterns of $D$ and $\widetilde{D}$. The detailed proof can be found in the Appendix.

THEOREM 1. *Let $D$ and $\widetilde{D}$ be a pair of complementary regular $2^{n-m}$ even designs as in (4). Then,*

$$
\begin{aligned}
A_{2u}(D) = {} & C_{2u} + C_{2u,0} + \gamma_{2u} + A_{2u}(\widetilde{D}) \\
& + \sum_{l=1}^{u-1}\left(\sum_{t=1}^{2u-1-2l} C_{2u,t+2l}\alpha_{k-1}(t) + \alpha_{k-1}(2(u-l))\right) \\
& \qquad \times \left(\binom{2^{k-1}-n}{2l} - A_{2l}(\widetilde{D})\right) \\
& + \sum_{l=2}^{u-1}\left(\sum_{t=1}^{2u-1-2l} C_{2u,t+2l}\gamma_{k-1}(t) + \gamma_{k-1}(2(u-l)) + C_{2u,2l}\right)A_{2l}(\widetilde{D})
\end{aligned}
$$

*for $u = 2, \ldots, [n/2]$, where $\alpha_{k-1}(s)$ and $\gamma_{k-1}(s)$ are defined in (5) and (6), respectively, $\gamma_{2u} = \sum_{s=3}^{2u} I_{[s \le 2^{k-1}-1]}\gamma_{k-1}(s)$, $C_i = 2^{-k}(P_i(0;n) - P_i(2^{k-1};n))$, $C_{ij} = (-1)^{i-[(i-j)/2]}\binom{n-2^{k-1}}{[(i-j)/2]}$ and $[x]$ is the largest integer less than or equal to $x$.*

By noting from (7) that

$$A_{2u}(D) = A_{2u}(\widetilde{D}) + \text{lower order terms},$$

we have the following result.

COROLLARY 1. *For $5N/16 < n < N/2$, a regular $2^{n-m}$ design has minimum aberration if and only if it is even and its complementary design in the maximal even design has minimum aberration among all $(N/2-n)$-factor even designs.*

This result was also obtained in Butler (2003) by comparing the design moments.

In the following section, we will use the identities in (7) to further investigate the structures of (weak) minimum aberration designs of resolution IV.

Xu and Cheng ([2008](#)) developed a general complementary design theory for doubling that can be applied to the case where $n \le 5N/16$. In this case, MA designs are projections of designs other than maximal even designs where, unlike the case $5N/16 < n < N/2$ studied here, wordlength patterns of the complementary designs alone are no longer sufficient to characterize minimum aberration projections.

**3. Weak minimum aberration $2^{n-m}$ designs of resolution IV.** In this section, we focus on weak minimum aberration $2^{n-m}$ designs with $5N/16 < n < N/2$. First, an explicit relationship between $A_4(D)$ and $A_4(\widetilde{D})$ is derived.

From Theorem [1](#), we have

$$
\begin{aligned}
(8) \qquad A_4(D) = {} & A_4(\widetilde{D}) + C_4 + C_{4,0} + \gamma_{k-1}(3) + \gamma_{k-1}(4) \\
& + (1 + \alpha_{k-1}(2)) \binom{2^{k-1} - n}{2}.
\end{aligned}
$$

The following corollary results from simplification of [(8)](#).

COROLLARY 2. *For each regular $2^{n-m}$ even design $D$,*

$$
(9) \qquad A_4(D) = A_4(\widetilde{D}) + \left( \binom{n}{4} - \binom{2^{k-1} - n}{4} \right) \Big/ (2^{k-1} - 3).
$$

From Corollary [2](#), we have the following result.

THEOREM 2. *For any regular $2^{n-m}$ even design $D$ of resolution* IV,

$$
(10) \qquad A_4(D) \ge \left( \binom{n}{4} - \binom{2^{k-1} - n}{4} \right) \Big/ (2^{k-1} - 3).
$$

If the complementary design $\widetilde{D}$ in [(4)](#) has resolution at least VI, then $D$ achieves the lower bound for $A_4(D)$ in Theorem [2](#), and is a weak minimum aberration design. Let $M_6(k)$ be the maximum number of factors that can be accommodated in an even fractional factorial design of resolution at least VI and run size $N = 2^k$. Then, for $N/2 - M_6(k) \le n \le N/2$, the lower bound in [(10)](#) can be attained, and a weak minimum aberration design can be constructed as the complement of an even fractional factorial design of resolution at least VI in the maximal even design. There is no explicit formula for $M_6(k)$, but some upper bounds on $M_6(k)$ can be derived from coding theory.

It is well known that the concepts of fractional factorial designs, wordlength pattern, resolution and defining relation have their counterparts in the context of linear codes. See MacWilliams and Sloane ([1977](#)) for basic concepts in algebraic coding theory. Let the defining words of a fractional factorial design be represented by binary row vectors. A regular $2^{n-m}$ fractional factorial design can be considered as an $[n, n-m]$ linear code, which is the null space of the $m \times n$ matrix whose rows are the $m$ independent defining words. Then, the defining relation of the design can be considered as its dual code, which is the $[n, m]$ linear code generated by the $m$ independent defining words. It follows, from the Varshamov lower bound in coding theory [Huffman and Pless (2003), page 87], that $M_6(k)$ must satisfy the following inequality:

$$(11) \qquad \binom{x-1}{1} + \binom{x-1}{2} + \binom{x-1}{3} + \binom{x-1}{4} \leq 2^k - 1.$$

For example, $M_6(6) = 7$, and, for $k = 6$, the maximum value of $x$ satisfying inequality (11) is also 7.

Let $[n, t, d]$ be a binary linear code with minimum distance $d$. Given a $[n, t, d]$ code with odd $d$, we can obtain an even $[n+1, t, d+1]$ code by appending a 0 to every codeword of even weight and a 1 to every codeword of odd weight. The new code is called the extended code. Obviously, extended codes are even. BCH codes are a family of multiple-error-correcting codes that were discovered by Bose, Ray-Chauhuri and Hocquenghen. From page 586 of MacWilliams and Sloane (1977), there exist BCH codes with parameters $[2^r, 2^r - 2r, 5]$ and $[2^r + 2^{[(r+1)/2]} - 1, 2^r + 2^{[(r+1)/2]} - 2r, 5]$. This implies the existence of the extended codes $[2^r + 1, 2^r - 2r, 6]$ and $[2^r + 2^{[(r+1)/2]}, 2^r + 2^{[(r+1)/2]} - 2r - 2, 6]$, which can be used to construct weak minimum aberration designs. We have the following theorem.

THEOREM 3. *Let $D^*$ be a weak minimum aberration $2^{n-m}$ design and $k = n - m$. If $k$ is odd and $2^{k-1} - 2^{(k-1)/2} - 1 \leq n \leq 2^{k-1}$, or $k$ is even and $2^{k-1} - 2^{(k-2)/2} - 2^{[k/4]} \leq n \leq 2^{k-1}$, then $A_4(D^*)$, the minimum number of words of length four, is*

$$(12) \qquad A_4(D^*) = \left( \binom{n}{4} - \binom{2^{k-1}-n}{4} \right) \Big/ (2^{k-1} - 3).$$

PROOF. Here, we only consider the case of even $k$. The case of odd $k$ can be handled similarly. Let $k = 2r + 2$ and $n^* = 2^r + 2^{[(r+1)/2]} = 2^{(k-2)/2} + 2^{[k/4]}$. Since the extended BCH code $[2^r + 2^{[(r+1)/2]}, 2^r + 2^{[(r+1)/2]} - 2r - 2, 6]$ exists, the maximum resolution of an even $2^{n^*-(n^*-k)}$ design is at least VI. Therefore, the maximum resolution of a complementary design $\tilde{D}$ in (4) is at least VI if the number of factors is greater than or equal to $2^{k-1} - 2^{(k-2)/2} - 2^{[k/4]}$. Equation (12) follows from Corollary 2.

A weak minimum aberration design attaining (12) can be obtained by first constructing the even design with an extended BCH code as its defining relation, and then taking its complement in the maximal even design.

From Theorem 3, for $k = 6$, equality (12) holds if $26 \leq n \leq 32$. Since $M_6(6) = 7$, (12) also holds for $n = 25$. However, even if $M_6(k)$ can be determined, the cases $5N/16 < n < N/2 - M_6(k)$ are not covered by Theorem 3. For $k = 6$, the maximum resolution of a complementary design $\widetilde{D}$ is IV when $20 < n < 25$. For those $n$ and $k$ that are not covered by Theorem 3, we will derive in the next section a lower bound on the minimum number of words of length four by using MacWilliams identities and linear programming techniques. Before doing that, we first examine the structure of the complementary design $\widetilde{D}$ when $D$ has minimum aberration.

Let $T$ and $\widetilde{T}$, respectively, be the factor representations of $D$ and $\widetilde{D}$ in (4). A word of length four in the defining relation of $D$ corresponds to a linearly dependent quadruple of points of $T$. From Corollary 2, $A_4(D)$ is minimized if and only if $A_4(\widetilde{D})$ is minimized; that is, $\widetilde{T}$ must contain the minimum number of linearly dependent quadruples. Let $M$ be an $m$-subset of $\mathrm{PG}(k-1, 2)$. The rank of $M$, denoted as $\mathrm{rank}(M)$, is defined as the maximum number of independent points of $M$. The following theorem indicates that, for a minimum aberration even design, the factor representation $\widetilde{T}$ of its complementary design must have maximum rank.   □

THEOREM 4.  *Let $M$ be an $m$-subset of $\mathrm{PG}(k-1, 2)$ containing the minimum number of linearly dependent quadruples among all $m$-subsets of $\mathrm{PG}(k-1, 2)$. Then, $M$ must have maximum rank.*

PROOF.  Let $M = \{\mathbf{a}_1, \ldots, \mathbf{a}_m\}$. If $\mathrm{rank}(M)$ is less than $k$, say, $\mathrm{rank}(M) = r < k$ [i.e., $M \subseteq \mathrm{PG}(r-1, 2)$], then $\mathrm{PG}(k-1, 2) \setminus \mathrm{PG}(r-1, 2)$ is not empty. Let $\{\mathbf{a}_1, \mathbf{a}_2, \mathbf{a}_3, \mathbf{a}_4\}$ be a linearly dependent quadruple of points in $M$. Since $\mathbf{a}_1 = \mathbf{a}_2 + \mathbf{a}_3 + \mathbf{a}_4$, $\mathrm{rank}(M \setminus \{\mathbf{a}_1\}) = r$. Let $\mathbf{a} \in \mathrm{PG}(k-1, 2) \setminus \mathrm{PG}(r-1, 2)$ and $M' = \{\mathbf{a}\} \cup M \setminus \{\mathbf{a}_1\}$. Obviously, $\mathrm{rank}(M') = r + 1 > r = \mathrm{rank}(M)$, and the number of linearly dependent quadruples in $M'$ is at least one less than that in $M$. This is a contradiction.   □

## 4. A lower bound on the minimum number of words of length four.
From Theorem 4, if a resolution IV even design $D$ has minimum aberration, then its complementary design $\widetilde{D}$ must have the maximum rank. To study minimum aberration designs for $5N/16 < n < N/2 - M_6(k)$, we assume that $D$ and $\widetilde{D}$ have rank $n - m$.

As discussed in Section 3, a regular $2^{n-m}$ design $D$ can be considered as an $[n, n-m]$ linear code, and the wordlength pattern of $D$, $\{A_i(D)\}$, is the same as the weight distribution of the dual code. Let $\{A'_i(D)\}$ be the weight distribution of $D$. MacWilliams identities in coding theory [MacWilliams and Sloane (1977)] provide a fundamental relationship between $\{A_i(D)\}$ and $\{A_i'(D)\}$,

$$(13) \qquad A_i(D) = 2^{m-n} \sum_{j=0}^{n} P_i(j;n) A_j'(D)$$

for $i = 0, \ldots, n$, where $P_i(j;n)$ is a Krawtchork polynomial.

Karpovsky (1979) derived another version of MacWilliams identities,

$$(14) \qquad A_i(D) = \frac{1}{i!}(S_i - C_i),$$

where $S_i = 2^{m-n} \sum_{j=0}^{n} (n-2j)^i A_j'(D)$ and $C_i$ is a constant.

For $i = 4$, $C_4 = n(3n-2)$. From (14), we have

$$(15) \qquad A_4(D) = \frac{1}{4!}\left(2^{m-n} \sum_{j=0}^{n} (n-2j)^4 A_j'(D) - n(3n-2)\right).$$

Since $D$ contains the vector $\mathbf{1} = (1, \ldots, 1)$, we have $A_0'(D) = A_n'(D)$ and $A_j'(D) = A_{n-j}'(D)$. For simplicity, let

$$(16) \qquad x_j = \begin{cases} A_j'(D), & \text{for } j = 1, \ldots, [(n-1)/2], \\ \frac{1}{2} A_j'(D), & \text{for } j = n/2, \text{ if } n \text{ is even.} \end{cases}$$

Replacing $A_j'(D)$ by (16), we can express equation (15) as

$$(17) \qquad A_4(D) = \frac{n^4}{12 * 2^{n-m}} - \frac{3n^2 - 2n}{24} + \frac{1}{12 * 2^{n-m}} \sum_{j=1}^{[n/2]} (n-2j)^4 x_j.$$

In (14), $C_2 = n$; thus, we have

$$(18) \qquad A_2(D) = \frac{1}{2!}\left(2^{m-n} \sum_{j=0}^{n} (n-2j)^2 A_j'(D) - n\right).$$

Also, it is not difficult to see that

$$(19) \qquad \sum_{j=0}^{n} A_j'(D) = 2^{n-m}.$$

Since $A_0'(D) = 1$ and $A_2(D) = 0$, the following equations result from replacing $A_j'(D)$ with (16) in (18) and (19), respectively:

$$(20) \qquad \sum_{j=1}^{[n/2]} x_j = 2^{n-m-1} - 1 \quad \text{and} \quad \sum_{j=1}^{[n/2]} (n-2j)^2 x_j = (2^{n-m-1} - n)n.$$

Thus, we can formulate a linear programming (LP) problem for bounding the minimum of $A_4(D)$. The LP problem is to find a vector $(x_1, \ldots, x_{[n/2]})$ that minimizes

(21)
$$f(x_1, \ldots, x_{[n/2]}) = \frac{n^4}{12 * 2^{n-m}} - \frac{3n^2 - 2n}{24}$$
$$+ \frac{1}{12 * 2^{n-m}} \sum_{j=1}^{[n/2]} (n-2j)^4 x_j,$$

subject to the linear constraints (20).

Based on the fact that the function $f$ in (21) assumes its minimum at an extreme point [see Gass (1985)], we derive the following theorem, whose proof can be found in the Appendix.

THEOREM 5. *When* $5N/16 < n < N/2 - M_6(k)$, *for any* $2^{n-m}$ *even design* $D$ *of resolution* IV,

(22)
$$A_4(D) \geq \mathrm{LB}(n, n-m)$$
$$= \frac{n^4}{12 * 2^{n-m}} - \frac{3n^2 - 2n}{24}$$
$$+ \frac{n^2(2^{n-m-1} - n)^2}{(2^{n-m-1} - 1) * 12 * 2^{n-m}}.$$

Similarly, a lower bound for $A_4(\widetilde{D})$ is $\mathrm{LB}(N/2 - n, n - m)$. From Theorem 5 and Corollary 2, we have the following theorem

THEOREM 6. *When* $5N/16 < n < N/2 - M_6(k)$, *for any* $2^{n-m}$ *even design* $D$ *of resolution* IV,

$$A_4(D) \geq \max\Big\{\mathrm{LB}(N/2 - n, n - m)$$
$$+ \Big(\binom{n}{4} - \binom{2^{k-1} - n}{4}\Big)\Big/(2^{k-1} - 3), \mathrm{LB}(n, n-m)\Big\}.$$

As discussed in Section 3, for $k = 6$, Theorem 3 does not cover the cases $20 < n < 25$. In these cases, lower bounds on $A_4(D)$ can be obtained from Theorem 6. Table 1 compares the lower bounds to the actual minimum values of $A_4$.

We can see that the lower bounds are very tight. For $k = 7$, Theorem 3 covers the cases $55 \leq n \leq 64$, and the lower bounds on $A_4(D)$ can be derived from Theorem 6 for $40 < n < 55$. The following table compares those bounds to the actual minimum values of $A_4$ obtained from Block and Mee (2005).

TABLE 1
*The minimum # of $A_4$ and their lower bounds*
*for 64-run designs*

| $n$ | Lower bound | min $A_4$ |
|---|---|---|
| 21 | 203 | 204 |
| 22 | 249 | 250 |
| 23 | 302 | 304 |
| 24 | 364 | 365 |

TABLE 2
*The minimum # of $A_4$ and their lower bounds*
*for 128-run designs*

| $n$ | Lower bound | min $A_4$ |
|---|---|---|
| 54 | 5181 | 5182 |
| 53 | 4795 | 4797 |
| 52 | 4431 | 4433 |
| 51 | 4089 | 4091 |
| 50 | 3766 | 3770 |
| 49 | 3463 | 3466 |
| 48 | 3179 | 3180 |
| 47 | 2912 | 2915 |
| 46 | 2662 | 2665 |
| 45 | 2428 | 2430 |
| 44 | 2210 | 2214 |
| 43 | 2007 | 2009 |
| 42 | 1818 | 1822 |
| 41 | 1643 | 1648 |

From Table 2, it appears that the bound is consistently close to the actual minimum. The biggest gap occurs at $n = 41$.

In principle, the techniques developed in this section can be generalized to other wordlengths. However, the complexity of the associated linear programming problems may vary with respect to different lengths.

**5. Complementary designs of some MA $2^{n-m}$ even designs.** Chen and Hedayat (1996) showed that, in the general setting (2), the complement of an MA design in a regular saturated design of resolution III must have minimum rank. Thus, in general, the complementary designs of MA designs in regular saturated designs of resolution III have the same structure as long as they have the same number of factors, regardless of the run size. In contrast, we have shown in Section 3 that the complement of an MA even design of resolution IV in the maximal even design must have maxi-mum rank. A consequence is that such complementary designs have different structures when the run sizes are different even if they have the same number of factors. This makes the cataloging of complementary designs of MA even designs more difficult. However, for small numbers of factors, we can still derive the structures of the complementary designs of some MA even designs.

For $n < N/2$, as pointed out in Section 2, the points of $T$ and $\widetilde{T}$ in (3) together can be considered to form a Euclidean geometry EG$(k-1, 2)$,

$$(23) \qquad \mathrm{EG}(k-1,2) = \{\underbrace{\mathbf{a}_1, \ldots, \mathbf{a}_n}_{T}, \underbrace{\mathbf{b}_1, \ldots, \mathbf{b}_{\widetilde{n}}}_{\widetilde{T}}\},$$

where $\widetilde{n} = N/2 - n$ is the number of factors of the complementary design determined by $\widetilde{T}$, and $\mathbf{b}_1, \ldots, \mathbf{b}_{\widetilde{n}}$ are the $\widetilde{n}$ points $\mathbf{a}_{n+N/2}, \ldots, \mathbf{a}_{N-1}$ in (3).

For $\widetilde{n} \leq k$, the factor representation $T$ corresponds to an MA design if $\widetilde{T}$ is such that

$$\widetilde{T} = \{\widetilde{n} \text{ independent points of } \mathrm{EG}(k-1, 2)\}.$$

In other words, an even design whose complement in the maximal even design contains $\widetilde{n}$ independent columns has minimum aberration.

For $\widetilde{n} = k+1$, when $\widetilde{n}$ is even, the complementary design of an MA design is a $2^{\widetilde{n}-1}$ design of resolution $\widetilde{n}$, and when $\widetilde{n}$ is odd, the complementary design of an MA design is a $2^{\widetilde{n}-1}$ design of resolution $\widetilde{n} - 1$.

For $\widetilde{n} = k+2$, without loss of generality, let the $k$ points $\{\mathbf{b}_1, \ldots, \mathbf{b}_k\}$ in (23) be independent points of EG$(k-1, 2)$, and $k = 3m + r$, where $0 \leq r < 3$. A factor representation $\widetilde{T}$ that corresponds to the complementary design of an MA design is

$$\widetilde{T} = \{\mathbf{b}_1, \ldots, \mathbf{b}_k, \mathbf{c}, \mathbf{d}\},$$

where, for $r = 0$, $\mathbf{c}$ and $\mathbf{d}$ are defined as

$$\mathbf{c} = \mathbf{b}_1 + \mathbf{b}_2 + \cdots + \mathbf{b}_{2m-1}, \qquad \mathbf{d} = \mathbf{b}_{m+1} + \mathbf{b}_{m+2} + \cdots + \mathbf{b}_{3m} + \mathbf{c};$$

for $r = 1$,

$$\mathbf{c} = \mathbf{b}_1 + \mathbf{b}_2 + \cdots + \mathbf{b}_{2m+1}, \qquad \mathbf{d} = \mathbf{b}_{m+1} + \mathbf{b}_{m+2} + \cdots + \mathbf{b}_{3m+1};$$

for $r = 2$,

$$\mathbf{c} = \mathbf{b}_1 + \mathbf{b}_2 + \cdots + \mathbf{b}_{2m+1}, \qquad \mathbf{d} = \mathbf{b}_{m+1} + \mathbf{b}_{m+2} + \cdots + \mathbf{b}_{3m+2} + \mathbf{c}.$$

For $\widetilde{n} = k+3$, following the notation in Chen and Wu (1991), let the defining relation of the complementary design of an MA design be

$$I = B_7 B_6 B_4 B_3 = B_7 B_5 B_4 B_2 = B_6 B_5 B_4 B_1.$$

Let $\tilde{n} = 7m + r$, $0 \leq r \leq 6$. These $B_i$ divide the $\tilde{n}$ letters into seven approximately equal blocks.

For $r = 0, 1$,

$$B_i = (im - m + 1)(im - m + 2) \cdots (im), \qquad i = 1, \ldots, 7.$$

For $r = 2$,

$$B_i = (im - m + 1)(im - m + 2) \cdots (im), \qquad i = 2, \ldots, 7,$$
$$B_1 = 1 \cdot 2 \cdots m(7m + 1)(7m + 2).$$

For $r = 3$,

$$B_i = (im - m + 1)(im - m + 2) \cdots (im)(7m + i), \qquad i = 1, 2,$$
$$B_j = (jm - m + 1)(jm - m + 2) \cdots (jm), \qquad j = 3, 4, 6, 7,$$
$$B_5 = (4m + 1)(4m + 2) \cdots 5m(7m + 3).$$

For $r = 4$,

$$B_i = (im - m + 1)(im - m + 2) \cdots (im)(7m + i), \qquad i = 1, 2, 3, 4,$$
$$B_j = (jm - m + 1)(jm - m + 2) \cdots (jm), \qquad j = 5, 6, 7.$$

For $r = 5$,

$$B_i = (im - m + 1)(im - m + 2) \cdots (im)(7m + i), \qquad i = 1, 2, 3,$$
$$B_j = (jm - m + 1)(jm - m + 2) \cdots (jm), \qquad j = 5, 6, 7,$$
$$B_4 = (3m + 1)(3m + 2) \cdots 4m(7m + 4)(7m + 5).$$

For $r = 6$,

$$B_i = (im - m + 1)(im - m + 2) \cdots (im)(7m + i), \qquad i = 1, 2, 3, 4,$$
$$B_j = (jm - m + 1)(jm - m + 2) \cdots (jm), \qquad j = 6, 7,$$
$$B_5 = (4m + 1)(4m + 2) \cdots 5m(7m + 5)(7m + 6).$$

For example, the complementary design of the $2^{8-2}$ design defined by $I = 1237 = 345678$ ($k = 6$) is an MA $2^{24-18}$ design. Similarly, the complementary design of the $2^{9-2}$ design defined by $I = 123458 = 345679$ ($k = 7$) is an MA $2^{55-48}$ design.

## APPENDIX: PROOFS

Note that, in (3), $\mathrm{PG}(k - 1, 2)$ is partitioned into three parts, one of which is itself a projective geometry. Such a partition also arises in the blocking of fractional factorial designs as studied in Chen and Cheng (1999). A key result there is useful for proving Theorem 1.

Following the notations in Chen and Cheng (1999), let $D_B(2^{n-m} : 2^r)$ be a $2^{n-m}$ design in $2^r$ blocks of size $2^{n-m-r}$ $(r < n - m)$. Then, $D_B(2^{n-m} : 2^r)$ can be viewed as a $2^{(n+r)-(m+r)}$ design, where the factors are divided into $n$ treatment factors and $r$ blocking factors. For a blocked design $D_B(2^{n-m} : 2^r)$, let $A_{i,0}(D_B)$ [resp., $A_{i,1}(D_B)$] be the number of treatment defining words (resp., block defining words) containing $i$ treatment letters. Since Chen and Cheng (1999) only considered designs in which none of the treatment main effects is aliased with other main effects or confounded with blocks, it was assumed that $A_{1,0} = A_{0,1} = A_{2,0} = A_{1,1} = 0$. The two vectors $W_t(D_B) = (A_{3,0}(D_B), A_{4,0}(D_B), \ldots)$ and $W_{bt}(D_B) = (A_{2,1}(D_B), A_{3,1}(D_B), \ldots)$ together are called the *split wordlength pattern* of $D_B$.

The set $\{D, D^F\}$ in (4) can be viewed as the design $D$ divided into $2^{k-1}$ blocks, denoted as $D_B(2^{n-m} : 2^{k-1})$. The set $\{\widetilde{D}, D^F\}$ in (4) represents the design $\widetilde{D}$ in $2^{k-1}$ blocks, denoted as $D_R$. The blocked design $D_R$ is called the *blocked residual design* of $D_B$. The wordlength pattern of $D$ in (4) corresponds to $W_t(D_B) = (A_{4,0}(D_B), 0, A_{6,0}(D_B), \ldots)$ of the blocked design $D_B$. Similarly, the wordlength pattern of its complementary design $\widetilde{D}$ is the same as $W_t(D_R) = (A_{4,0}(D_R), 0, A_{6,0}(D_R), \ldots)$ of the blocked residual design $D_R$. Chen and Cheng (1999) showed that, in general, $A_{i,0}(D_B)$ can be written in terms of the split wordlength pattern of $D_R$. The following result is Theorem 1 of Chen and Cheng (1999).

THEOREM 7. *Let $\{A_{i_1,b}(D_B)\}$ and $\{A_{i_1,b}(D_R)\}$ be the split wordlength patterns of a blocked design $D_B(2^{n-m} : 2^r)$ and its blocked residual design $D_R$, respectively. Then,*

$$
\begin{aligned}
A_{i,0}(D_B) = {} & C_i + C_{i0} + \sum_{s=3}^{i} I_{[s \leq 2^r - 1]}\gamma_r(s) \\
& + \sum_{s=2}^{i-1}\bigg(\sum_{t=1}^{i-1-s} C_{i,t+s}I_{[t \leq 2^r - 2]}\alpha_r(t) \\
& \qquad\qquad + (-1)^i I_{[i-s \leq 2^r - 2]}\alpha_r(i-s)\bigg)A_{s,1}(D_R) \\
& + \sum_{s=2}^{i-1}\bigg(\sum_{t=1}^{i-1-s} C_{i,t+s}I_{[t \leq 2^r - 2]}\gamma_r(t) \\
& \qquad\qquad + (-1)^i I_{[i-s \leq 2^r - 2]}\gamma_r(i-s)\bigg)A_{s,0}(D_R) \\
& + \sum_{s=3}^{i} C_{i,s}(A_{s,0}(D_R) + I_{[s>2^r]}A_{s-2^r+1,0}(D_R))
\end{aligned}
$$
(24)

*for $i = 3, \ldots, n$, where $C_{ij} = (-1)^{i-[(i-j)/2]}\binom{n-2^{k-1}}{[(i-j)/2]}$, $C_i = 2^{-k}[P_i(0;n) - P_i(2^{k-1};n)]$, $[x]$ is the largest integer less than or equal to $x$ and $I_{[\cdot]}$ is the indicator function that takes the value 1 or 0 depending on whether condition $[\cdot]$ is true or not.*

Note that the term $\sum_{s=3}^{i} I_{[s \leq 2^r - 1]} \gamma_r(s)$ is missing from Theorem 1 of Chen and Cheng ([1999](#)), which is an error.

PROOF OF THEOREM [1]. In ([24](#)), $I_{[t \leq 2^r - 2]} = I_{[i-s \leq 2^r-2]} = 1$ and $I_{[s > 2^r]} = 0$ when $r = k-1$. The blocked residual design $D_R$ determined by $\{\widetilde{D}, D^F\}$ in ([4](#)) has $A_{3,0}(D_R) = 0$, $A_{2,1}(D_R) = \binom{2^{k-1}-n}{2}$, $A_{3,1}(D_R) = 0$, $A_{4,0}(D_R) + A_{4,1}(D_R) = \binom{2^{k-1}-n}{4}$, $A_{5,0}(D_R) = 0$, $A_{5,1}(D_R) = 0$, etc. Therefore,

$$(25) \qquad A_{s,1}(D_R) = \begin{cases} \binom{2^{k-1}-n}{s} - A_{s,0}(D_R), & \text{if } s \text{ is even,} \\ 0, & \text{if } s \text{ is odd.} \end{cases}$$

Since $D$ in ([4](#)) is an even design, the subscript of $A_i(D)$ can be denoted as $i = 2u$, where $u = 2, \ldots, [n/2]$. Let $\gamma_{2u} = \sum_{s=3}^{2u} I_{[s \leq 2^{k-1}-1]} \gamma_{k-1}(s)$. Similarly, let $s = 2l$ where $l = 1, 2, \ldots, u-1$ and replace $\{A_{2l,1}(D_R)\}$ in ([24](#)) by ([25](#)); then, we obtain equation ([7](#)). □

PROOF OF THEOREM [5]. Let

$$(26) \qquad h(x_1, \ldots, x_{[n/2]}) = \sum_{j=1}^{[n/2]} (n-2j)^4 x_j.$$

To minimize ([21](#)) is equivalent to finding the minimum of ([26](#)). From Chapter 3 of Gass ([1985](#)), $h(x_1, \ldots, x_{[n/2]})$ assumes its minimum at an extreme point $\mathbf{x} = (0, \ldots, 0, x_l, 0, \ldots, 0, x_g, 0, \ldots, 0)$, where $x_l$ and $x_g$ are nonnegative and are subject to the linear constraints ([20](#)); that is,

$$(27) \qquad \begin{aligned} x_l + x_g &= 2^{n-m-1} - 1 \quad \text{and} \\ (n-2l)^2 x_l + (n-2g)^2 x_g &= (2^{n-m-1}-n)n. \end{aligned}$$

The solution to ([27](#)) is

$$(28) \qquad \begin{aligned} x_l &= \frac{(2^{n-m-1}-1)(n-2g)^2 - (2^{n-m-1}-n)n}{(n-2g)^2 - (n-2l)^2}, \\ x_g &= \frac{(2^{n-m-1}-1)(n-2l)^2 - (2^{n-m-1}-n)n}{(n-2g)^2 - (n-2l)^2}. \end{aligned}$$

Plugging (28) into (26), we have

$$
\begin{aligned}
h(x_1, &\ldots, x_{[n/2]}) \\
&= h(l, g) \\
&= \frac{(n-2l)^4[(2^{n-m-1}-1)(n-2g)^2 - (2^{n-m-1}-n)n]}{(n-2g)^2 - (n-2l)^2} \\
&\quad + \frac{(n-2g)^4[(2^{n-m-1}-n)n - (2^{n-m-1}-1)(n-2l)^2]}{(n-2g)^2 - (n-2l)^2}.
\end{aligned}
$$
(29)

The function (29) assumes its minimum at the point $(l^*, g^*)$ that satisfies the equations

$$
\frac{\partial h(l, g)}{\partial l} = 0,
$$

$$
\frac{\partial h(l, g)}{\partial g} = 0.
$$

The minimum of function (29) is achieved at

$$
l^* = 1 + \frac{1}{2}\left(n - \sqrt{\frac{(2^{n-m-1}-n)n}{(2^{n-m-1}-1)}}\right),
$$

$$
g^* = \frac{1}{2}\left(n - \sqrt{\frac{(2^{n-m-1}-n)n}{(2^{n-m-1}-1)}}\right).
$$

Thus, the minimum value of (29) is $(2^{n-m-1}-n)n^2/(2^{n-m-1}-1)$. Replacing the minimum value in (21), we obtain the inequality (22). □

DIVISION OF BIOSTATISTICS
  AND BIOINFORMATICS
SCHOOL OF MEDICINE
UNIVERSITY OF MARYLAND
660 WEST REDWOOD STREET
BALTIMORE, MARYLAND 21201
USA
E-MAIL: hchen@epi.umaryland.edu

DEPARTMENT OF STATISTICS
UNIVERSITY OF CALIFORNIA, BERKELEY
BERKELEY, CALIFORNIA 94720-3860
USA
E-MAIL: cheng@stat.berkeley.edu